\documentclass[final,12pt]{elsarticle}
\usepackage{srcltx}
\usepackage{eurosym}
\usepackage{mathtools}
\usepackage{amsmath}
\usepackage{amsfonts}
\usepackage{amssymb}
\usepackage{amsthm}
\usepackage{graphicx}
\usepackage{mathrsfs}
\usepackage{xcolor}
\usepackage{exscale}
\usepackage{latexsym}

\usepackage[colorlinks,plainpages=true,pdfpagelabels,hypertexnames=true,colorlinks=true,pdfstartview=FitV,linkcolor=blue,citecolor=red,urlcolor=black]{hyperref}
\PassOptionsToPackage{unicode}{hyperref}
\PassOptionsToPackage{naturalnames}{hyperref}
\usepackage{enumerate}
\usepackage[shortlabels]{enumitem}
\usepackage{bookmark}
\usepackage{wasysym}
\usepackage{esint}
\usepackage[ddmmyyyy]{datetime}
\usepackage[margin=2cm]{geometry}
\makeatletter
\g@addto@macro\normalsize{%
	\setlength\abovedisplayskip{4pt}
	\setlength\belowdisplayskip{4pt}
	\setlength\abovedisplayshortskip{4pt}
	\setlength\belowdisplayshortskip{4pt}
}
\numberwithin{equation}{section}
\everymath{\displaystyle}
\usepackage[capitalize,nameinlink]{cleveref}
\crefname{section}{Section}{Sections}
\crefname{subsection}{Subsection}{Subsections}
\crefname{condition}{Condition}{Conditions}
\crefname{hypothesis}{Hypothesis}{Conditions}
\crefname{assumption}{Assumption}{Assumptions}
\crefname{lemma}{Lemma}{Lemmas}
\crefname{definition}{Definition}{Definitions}

\crefformat{equation}{\textup{#2(#1)#3}}
\crefrangeformat{equation}{\textup{#3(#1)#4--#5(#2)#6}}
\crefmultiformat{equation}{\textup{#2(#1)#3}}{ and \textup{#2(#1)#3}}
{, \textup{#2(#1)#3}}{, and \textup{#2(#1)#3}}
\crefrangemultiformat{equation}{\textup{#3(#1)#4--#5(#2)#6}}%
{ and \textup{#3(#1)#4--#5(#2)#6}}{, \textup{#3(#1)#4--#5(#2)#6}}%
{, and \textup{#3(#1)#4--#5(#2)#6}}

\Crefformat{equation}{#2Equation~\textup{(#1)}#3}
\Crefrangeformat{equation}{Equations~\textup{#3(#1)#4--#5(#2)#6}}
\Crefmultiformat{equation}{Equations~\textup{#2(#1)#3}}{ and \textup{#2(#1)#3}}
{, \textup{#2(#1)#3}}{, and \textup{#2(#1)#3}}
\Crefrangemultiformat{equation}{Equations~\textup{#3(#1)#4--#5(#2)#6}}%
{ and \textup{#3(#1)#4--#5(#2)#6}}{, \textup{#3(#1)#4--#5(#2)#6}}%
{, and \textup{#3(#1)#4--#5(#2)#6}}

\crefdefaultlabelformat{#2\textup{#1}#3}
\numberwithin{equation}{section}

\newtheorem{theorem} {Theorem}[section]

\newtheorem{counter example}{Counter Example}[section]
\newtheorem{remark} {Remark}[section]
\newtheorem{definition} {Definition}[section]

\newtheorem{claim} {Claim}[section]

\def\CC{{\rm \kern.24em \vrule width.02em height1.4ex depth-.05ex \kern-.26emC}}

\def\TagOnRight

\def\AA{{it I} \hskip-3pt{\tt A}}

\def\QQ{\rlap {\raise 0.4ex \hbox{$\scriptscriptstyle |$}} {\hskip -0.1em Q}}

\makeatletter
\newcommand{\vo}{\vec{o}\@ifnextchar{^}{\,}{}}
\makeatother

\def\YYint#1#2#3{{\setbox0=\hbox{$#1{#2#3}{\iint}$}
		\vcenter{\hbox{$#2#3$}}\kern-.50\wd0}}

\def\XXint#1#2#3{{\setbox0=\hbox{$#1{#2#3}{\int}$}
		\vcenter{\hbox{$#2#3$}}\kern-.50\wd0}}

\makeatletter
\def\namedlabel#1#2{\begingroup
	\def\@currentlabel{#2}%
	\label{#1}\endgroup
}
\makeatother
\makeatletter
\newcommand{\rmh}[1]{\mathpalette{\raisem@th{#1}}}
\newcommand{\raisem@th}[3]{\hspace*{-1pt}\raisebox{#1}{$#2#3$}}
\makeatother

\newcommand{\descitem}[2]{\item[#1] \label{#2}}
\newcommand{\descref}[2]{\hyperref[#1]{\textnormal{\textcolor{black}{(}\textcolor{blue}{\bf #2}\textcolor{black}{)}}}}

\newcommand{\dref}[2]{\hyperref[#1]{\textcolor{black}{(}\textcolor{blue}{\bf #2}\textcolor{black}{)}}}
\newcommand{\be} {\begin{eqnarray}}
\newcommand{\ee} {\end{eqnarray}}
\newcommand{\Bea} {\begin{eqnarray*}}
	\newcommand{\Eea} {\end{eqnarray*}}
\newcommand{\pa} {\partial}
\newcommand{\re}{\mathbb{R}}

\newcommand{\al} {\alpha}
\newcommand{\rr}{\rightarrow}

\newcommand{\dip}{\displaystyle}

\newcommand{\de} {\delta}
\newcommand{\g} {\gamma}

\newcommand{\p}  {\prime}

\newcommand{\la} {\lambda}
\newcommand{\si} {\sigma}

\newcommand{\f}{\infty}

\newcommand{\lab} {\label}

\newcounter{whitney}
\refstepcounter{whitney}

\newcounter{ineqcounter}
\refstepcounter{ineqcounter}
\makeatletter
\def\ps@pprintTitle{%
	\let\@oddhead\@empty
	\let\@evenhead\@empty
	\def\@oddfoot{}%
	\let\@evenfoot\@oddfoot}
\makeatother
\usepackage[doublespacing]{setspace}
\usepackage[titletoc,toc,page]{appendix}

\makeatletter
\newcommand{\refcheckize}[1]{%
	\expandafter\let\csname @@\string#1\endcsname#1%
	\expandafter\DeclareRobustCommand\csname relax\string#1\endcsname[1]{%
		\csname @@\string#1\endcsname{##1}\wrtusdrf{##1}}%
	\expandafter\let\expandafter#1\csname relax\string#1\endcsname
}
\makeatother

\refcheckize{\cref}
\refcheckize{\Cref}


\makeatletter
\newcommand{\mainsectionstyle}{%
	\renewcommand{\@secnumfont}{\bfseries}
	\renewcommand\section{\@startsection{section}{2}%
		\z@{.5\linespacing\@plus.7\linespacing}{-.5em}%
		{\normalfont\bfseries}}%
}
\makeatother
\usepackage{pgf,tikz}
\usetikzlibrary{arrows}
\usetikzlibrary{decorations.pathreplacing}

\usepackage{xpatch}
\xpatchcmd{\MaketitleBox}{\hrule}{}{}{}
\xpatchcmd{\MaketitleBox}{\hrule}{}{}{}

\date{}
\begin{document}
\begin{frontmatter}
	
	\title{Optimal jump set in hyperbolic conservation laws}
	
	\author[myaddress]{Shyam Sundar Ghoshal\tnoteref{thankssecondauthor}}
	\ead{ghoshal@tifrbng.res.in}
	\tnotetext[thankssecondauthor]{Supported in part by Inspire Research Grant.}

	\author[myaddress]{Animesh Jana}
	\ead{animesh@tifrbng.res.in }
	
	\address[myaddress]{Centre for Applicable Mathematics,Tata Institute of Fundamental Research, Post Bag No 6503, Sharadanagar, Bangalore - 560065, India.}
	\begin{abstract}
		This paper deals with some qualitative properties of entropy solutions to hyperbolic conservation laws. In \cite{DOW} the jump set of entropy solution to conservation laws has been introduced. We find an entropy solution to scalar conservation laws  for which the jump set is not closed, in particular, it is dense in a space-time domain. 
		 In the later part of this article we obtain a similar result for the hyperbolic system. 
 We give two different approaches for scalar conservation laws and hyperbolic system to obtain the results. For the scalar case, obtained solutions are more explicitly calculated.   
	\end{abstract}

	\end{frontmatter}
\tableofcontents
\section{Introduction}

The aim of this article has two folds. First, we obtain entropy solutions for scalar conservation laws, which have discontinuity on a dense set. Then we extend the result for strictly hyperbolic system (see definition \ref{SHS}). 
 First part of this section is devoted for the discussion about the problem and state of the art for scalar conservation laws.  Later in subsection \ref{intro-system}, we discuss about the relevant results in the existing literature for hyperbolic system. 

\subsection{Discussion on scalar conservation laws}
We consider the following multi-dimensional scalar conservation laws
\begin{eqnarray}
\frac{\pa}{\pa t}u+\sum\limits_{j=1}^{d}\frac{\pa}{\pa x_j}f_j(u)&=&0\ \ \ \mbox{ for } (x,t)\in\re^d\times\re_+,\lab{I1}\\
u(x,0)&=&u_0\ \ \mbox{for }x\in\re^d, \lab{I2}
\end{eqnarray}
where $u:\re^d\times\re_+\rr\re$ and $f=(f_1,\cdots,f_d)\in C^2\left(\re,\re^d\right)$ for $d\geq1$. It is well known that discontinuity may appear in the solution even for smooth initial data. Our main interest is to show the optimality of the discontinuity set of entropy solution. $u\in L^{\f}$ is called an entropy solution to (\ref{I1}) if it is a weak solution satisfying the Kru\v{z}kov \cite{Kruzkov} entropy condition in the sense of distribution:
\begin{equation}
\frac{\pa}{\pa t}|u-k|+\sum\limits_{j=1}^{d}\frac{\pa}{\pa x_j}\left(sgn(u-k)(f_j(u)-f_j(k))\right)\leq 0 \ \mbox{ for each }k\in\re.
\end{equation}
Suppose $u\in BV(\Omega)$ is an entropy solution to (\ref{I1}) where $\Omega\subset\re^d\times\re_+$ is an open set. It is well known (see for instance \cite{Ambrosio2}) that the \textit{approximate jump set} (see definition \ref{approx-jump-set}) of $u$ is $\mathcal{H}^{d}$--\textit{rectifiable} (see definition \ref{rectifiable}). It was predicted (see page 24, \cite{Dafermos1}) that the entropy solution can be discontinuous over a dense set. In the present article, we settle this question for any $C^2$ flux $f$.
\par In 1994, Lions, Perthame and Tadmor \cite{LPT} gave the notion of kinetic formulation for scalar conservation laws. Kinetic formulation for (\ref{I1}), (\ref{I2}) typically look like
\begin{eqnarray}
\frac{\pa}{\pa t}g(t,x,\xi)+\sum\limits_{k=1}^{d}a_{k}(\xi)\frac{\pa}{\pa x_k}g(t,x,\xi)&=&\frac{\pa}{\pa \xi}\mu(t,x,\xi)\ \ \mbox{in }\mathcal{D}^{\p}(\re_+\times\re^d\times\re),\lab{kinetic}\\
g(0,x,\xi)&=&\chi(\xi,u_0(x))\ \ \mbox{ for }(x,\xi)\in\re^d\times\re,
\end{eqnarray}
where $(a_1(\xi),\cdots,a_d(\xi))=f^{\p}(\xi)$ and the chi-function $\chi(\xi;u)$  is defined as follows:
\begin{equation}
\chi(\xi;u):=\left\{\begin{array}{lll}
&1&\mbox{if }0<\xi<u,\\
&-1&\mbox{if }u<\xi<0,\\
&0&\mbox{otherwise.}
\end{array}\right.
\end{equation}
In the kinetic equation (\ref{kinetic}), $\mu$ is a non-negative bounded measure which is usually referred as \textit{defect measure} in the literature. It is easy to see that  a smooth solution to (\ref{I1}) satisfies (\ref{kinetic}) with $0$ in the right hand side. Note that the measure $\mu$ captures the discontinuity of the solution. As we have discussed in the previous paragraph,  the approximate jump set of a $BV$ entropy solution is rectifiable. In general, the entropy solution  to scalar conservation laws need not be in $BV$ space (see for instance \cite{GJ}) unless the initial data is in $BV$ or the flux is uniformly convex in the $d=1$ case.  In \cite{DOW}, De Lellis, Otto and Westdickenberg derived a BV like structure of entropy solution for non degenerate flux $f$. In that article authors defined the set $J$ as follows
\begin{equation}\lab{jump}
J:=\left\{y\in\re^d\times\re_+;\limsup\limits_{r\downarrow0}\frac{\|\mu\|(B_r(y)\times\re)}{r^{d}}>0\right\},
\end{equation} 
where $\|\mu\|$ denotes the total variation of the defect measure $\mu$. Note that $J$ can be thought as positive upper $\mathcal{H}^{d}$ density of $\nu$ where $\nu(A):=\|\mu\|(A\times\re)$ for all Borel set $A\subset\re^{d}\times\re_+$. The authors have shown in \cite{DOW} that $J$ is a rectifiable set in the sense that  it is contained in a union of countable Lipschitz graphs. In \cite{Otto}, it has been shown that the set non-Lebesgue points has Hausdorff dimension atmost one for a larger class of solutions to Burger's equation.  The jump set $J$ of the entropy solution constructed in the present article is subset of countable union of disjoint hyperplanes. The union of those hyperplanes form a dense subset of $\re^d\times\re_+$.

\par For the rest of our discussion in this subsection we need to recall a non-degenerate flux condition from \cite{LPT}. We say the flux is \textit{non-degenerate} if there exist $\al\in(0,1]$ and a constant $C\geq0$ such that the following holds for $a:=f^{\p}$,
\begin{equation}\lab{flux-condition}
\sup\limits_{\tau^2+|\xi|^2=1} \mathcal{L}^1\left(\{|v|<R_0,\ |\tau +a(v)\cdot \xi|<\delta\}\right)<C \delta^\al\ \mbox{ for any }\de>0.
\end{equation}
Here $\mathcal{L}^1$ denotes the one dimensional Lebesgue measure. We use the notation $cl(A)$ as the closure of a set $A\subset\re^{N}$ in the standard topology for some $N\geq1$.

\par The aim of this  present article is related to the first  open problem mentioned in \cite[page 44]{Perthame}. One of the old  conjectures concerns the concentration of  $\mu$ on the jump set $J$. 
In a recent article \cite{Silvestre}, Silvestre proved the following theorem.
\begin{theorem}[Silvestre, \cite{Silvestre}]
Let $f$ be satisfying the non-degeneracy flux condition (\ref{flux-condition}). Let $u$ be an entropy solution to (\ref{kinetic}). Let $\mu$ be its kinetic dissipation measure and $J$ be its jump set. Then 
\begin{equation*}
\mu((\Omega\setminus{cl(J)})\times\re)=0.
\end{equation*}
\end{theorem}
 Now it remains to answer the following open question.\textbf{ Does $J$ differ from its closure $cl(J)$?} (see \cite[page 2]{Silvestre}). {\it In this paper, we find an entropy solution such that the corresponding jump set $J$ is not closed (see Theorem \ref{theorem1}). Together  with Silvestre's result and  Theorem \ref{theorem1} completely answers the conjecture  concerning the jump set of conservation laws.}
 Our result also tells that Silvestre's \cite{Silvestre} results cannot be improved that is to say that  there is no entropy dissipation outside $cl(J)$.
 
  Schaeffer \cite{Schaeffer} showed that $J$ is a closed set (in fact, union of finitely many curves) for a class of rapidly decaying functions. In \cite{AGG-structure} first author with Adimurthi and Gowda showed that $J$ can be union of infinitely many curves for $C_c^\f$ initial data. They have shown the existence of infinitely many \textit{asymptotically single shock packets}(see definition in \cite{AGG-structure}) for convex flux. These \textit{asymptotically single shock packets} form an open set. Hence they cannot be dense. For more detail check \cite{AGG-structure}. Dafermos \cite{Dafermos} proved that for $C^\f$ data and non-convex flux the entropy solution to (\ref{I1}) is piecewise $C^\f$. In this paper, we find an entropy solution $u$ for any $C^2$ flux such that the jump set $J$ is not closed. Note that by our method we can find a solution to scalar conservation laws with dense discontinuities. To the best of our knowledge this is the first time such a result (Theorem \ref{theorem1}) has been obtained even for Burger's equation. 
  
\subsection{Discussion on system of conservation laws}\lab{intro-system}
In this subsection, we mainly focus on the following hyperbolic system of conservation laws in one dimension.
\begin{equation}\lab{HS1}
\frac{\pa}{\pa t}U+\frac{\pa}{\pa x}F(U)=0 \ \mbox{ for } (x,t)\in\re\times\re_+.
\end{equation}
where $U:\re\times\re_+\rr\re^n$.

 Lax \cite{Lax} gave the existence of an admissible solution to Riemann problem for one dimensional hyperbolic system of conservation laws. The existential results for an initial data with small total are available due to Glimm's scheme \cite{Glimm}, front tracking algorithm by Bressan \cite{Bressan} and vanishing viscosity method by Bianchini and Bressan \cite{BB}. These solutions are admissible in the sense of Lax. Uniqueness result holds among the solutions obtained by the above mentioned procedures (see \cite{BF,Bressan3}).  Note that the general uniqueness result is still unavailable. Stability of entropy solutions has been shown for some special classes by Diperna \cite{Diperna}.  In the present paper, we construct solutions by an approximation through the  only finitely Riemann problem. These solutions satisfy the following Lax's admissible  condition. 
\begin{equation}
 \la_i(u_+)\leq\la\leq\la_i(u_-). 
\end{equation}

Previously, in \cite{Bressan2} the authors considered a $2\times2$ hyperbolic system with the following assumption
\begin{equation}\lab{BC1}
-R<\la_1(u)<-r<0<r<\la_2(u)<R, \mbox{ for all }u\in\Omega,
\end{equation} 
for some $r<R\in\re$ and open set $\Omega\subset\re^n$. They gave an initial data such that corresponding entropy weak solutions have small total variation and set of shocks is dense. Therefore it can not reach a constant state in finite time by any kind of boundary control. Note that this construction \cite{Bressan2} is specific with the assumption (\ref{BC1}). One can not do the same for scalar case and general strictly hyperbolic system. In this article, we consider strictly hyperbolic system with the assumption that there is atleast one $i\in\{1,\cdots,n\}$ such that the $i$-th characteristic field is either genuinely non-linear or linearly degenerate. Then we construct an initial data such that the shocks of admissible solution are dense.

\par  Lions, Perthame and Tadmor \cite{LPT1} considered the following $2\times2$ hyperbolic system of conservation laws
\begin{eqnarray}
\frac{\pa}{\pa t}\varrho+\frac{\pa}{\pa x}(\varrho u)&=&0,\lab{E1}\\
\frac{\pa}{\pa t}(\varrho u)+\frac{\pa}{\pa x}(\varrho u^2+\kappa\varrho^{\g})&=&0\lab{E2},
\end{eqnarray}
where $\g>1$ and $\kappa=\frac{(\g-1)^2}{4\g}$. In that article they gave the following kinetic formulation of the equations((\ref{E1}) and (\ref{E2})) of gas dynamics, 
\begin{equation}
\frac{\pa}{\pa t}\Phi+\frac{\pa}{\pa x}\left[(\theta\xi+(1-\theta)u)\Phi\right]=\frac{\pa^2}{\pa \xi^2}m(t,x,\xi).\lab{kinetic-gas}
\end{equation}
where $\theta=\frac{\g-1}{2}$ and the function $\Phi$ is defined as follows,
\begin{equation}
\Phi(\varrho;w):=\left(\varrho^{\g-1}-w^2\right)_+^\la,
\end{equation}
for $\la=\frac{3-\g}{2(\g-1)}$. One can compare the function $\Phi$ with the chi-function $\chi$ defined for the scalar case.  For the system (\ref{E1})-(\ref{E2}) a smooth solution satisfies (\ref{kinetic-gas}) with right hand side $0$ (see \cite{LPT1}). This indicates that the support of measure $m$ relates the discontinuity set of the corresponding entropy solution to the system (\ref{E1})-(\ref{E2}). Another open problem in (\cite{Perthame}, page 44) is the characterization of $m$. Theorem \ref{theorem2} can be applied for the system (\ref{E1})-(\ref{E2}). The entropy solution obtained in Theorem \ref{theorem2} is in $BV(\re\times[0,T])$ for some $0<T<1$ and it has discontinuities on a dense set. From the theory of $BV$ functions we know that the approximate jump set is $\mathcal{H}^1$ rectifiable. For the solution $U$ in Theorem \ref{theorem2}, the defect measure $m$ is a Radon measure supported on a union $\mathcal{B}$ of countably many line segments and $\mathcal{B}\neq cl(\mathcal{B})$. 

\par The article is organized as follows. In section \ref{prelim} we give some preliminary definitions and results. In section \ref{main-results} we state our main theorems. We discuss about the method and ideas in section \ref{method}. We put some preliminary results on $BV$ functions in subsection \ref{prelim}. We first prove the one dimensional case of Theorem \ref{theorem1} in subsection \ref{scalar-1d} and then multi dimensional case in subsection \ref{scalar-multi}. Then we prove Theorem \ref{theorem2} in section \ref{construction-system}.

\section{Preliminaries and notations}\lab{prelim}
Here we give some preliminary results and definition related to $BV$ functions. Throughout the section $\Omega$ is an open subset of $\re^d$.
\begin{definition}\lab{bv}
	Let $u\in L^1(\Omega)$. We say $u\in BV(\Omega)$ if there exists a finite signed Borel measure $\nu_i$ for each $i\in\{1,\cdots,d\}$ such that
	\begin{equation*}
	\int\limits_{\Omega}u\frac{\partial \psi}{\partial x_i}{d}x=-\int\limits_{\Omega}\psi{d}\nu_i,
	\end{equation*}
	holds for all $\psi\in C_c^{\f}(\Omega)$.
\end{definition}
Next, we recall the definitions of \textit{rectifiable sets, approximate limit} and \textit{approximate jump set} from \cite{Ambrosio2}.
\begin{definition}[\textbf{Rectifiable sets}]\lab{rectifiable}
	Let $A\subset\re^d$ be an $\mathcal{H}^k$ measurable set for $k\in\{0,1,\cdots,d\}$. $A$ is called \textit{countably $\mathcal{H}^k$--rectifiable} if there exist countably many Lipschitz maps $\phi_j:\re^k\rr\re^d$ such that
	\begin{equation*}
	\mathcal{H}^k\left(A\setminus\bigcup\limits_{j=1}^{\f}\phi_j(\re^k)\right)=0.
	\end{equation*}
	The set $A$ is called \textit{$\mathcal{H}^k$--rectifiable} if $A$ is {countably $\mathcal{H}^k$--rectifiable} and $\mathcal{H}^k(A)<\f$.
\end{definition}
\begin{definition}[\textbf{Approximate limit}]
	Let $v\in [L^1_{loc}(\Omega)]^m$. Let $x\in\Omega$. $v$ is said to have approximate limit at $x$ if there exists a $v_0\in\re^m$ such that 
	\begin{equation}
	\lim\limits_{r\downarrow0}\frac{1}{\mathcal{L}^d(B_r(x))}\int\limits_{B_r(x)}|v(y)-v_0|dy=0.
	\end{equation}
\end{definition}
Note that if $v$ has an \textit{approximate limit} at $x\in\Omega$ then $x$ is a Lebesgue point.
\begin{definition}[\textbf{Approximate discontinuity set}]\lab{discont}
	Let $u\in [L^1_{loc}(\Omega)]^m$. Let $S_u$ be the set of all points $x\in\Omega$ such that
	\begin{equation*}
	\lim\limits_{r\rr0}\frac{1}{\mathcal{L}^d(B_r(x))}\int\limits_{B_r(x)}|u(y)-z|dy>0,
	\end{equation*} 
	for all $z\in\re$. Note that if $x\notin S_u$ then $u$ has an \textit{approximate limit} at $x$. The set $S_u$ will be referred as approximate discontinuity set of $u$.	
\end{definition}

\begin{definition}[\textbf{Approximate jump set}]\lab{approx-jump-set}
	Any point $x\in\Omega$ is called \textit{approximate jump point} if there exists $a,b\in\re$ and $\nu\in S^{d-1}$ such that $a\neq b$ and the following holds
	\begin{equation}\lab{111}
	\lim\limits_{r\rr0}\frac{1}{\mathcal{L}^d(B^{-}_r(x,\nu))}\int\limits_{B^{-}_r(x,\nu)}|u(y)-a|dy=0\mbox{ and }\lim\limits_{r\rr0}\frac{1}{\mathcal{L}^d(B^{+}_r(x,\nu))}\int\limits_{B^{+}_r(x,\nu)}|u(y)-b|dy=0,
	\end{equation}
	where $B^{+}_r(x,\nu)=\{y\in B_r(x);(y-x)\cdot\nu>0\},B^{-}_r(x,\nu)=\{y\in B_r(x);(y-x)\cdot\nu<0\}$ and $S^{d-1}=\{\xi\in\re^d;|\xi|=1\}$. $a,b,\nu$ are uniquely determined by (\ref{111}) up to a permutation of $(a,b)$ and a change of sign of $\nu$.  We denote the triplet $(a,b,\nu)$ by $(u^+(x),u^-(x),\nu(x))$. The set of approximate jump points is denoted by ${J}_u$.
\end{definition}
\begin{theorem}[Federer\cite{Federer}--Vol'pert\cite{Volpert}]
	Let $u\in [L^1_{loc}(\Omega)]^m$. Let $S_u,J_u$ be as in definition \ref{discont} and definition \ref{approx-jump-set} respectively. Then $S_u$ is $\mathcal{H}^{d-1}$--rectifiable and $\mathcal{H}^{d-1}(S_u\setminus J_u)=0$. Moreover, we get
	\begin{equation*}
	Du\left|_{J_u}\right.=(u^+-u^-)\otimes\nu_u\mathcal{H}^{d-1}\left|_{J_u}\right.,
	\end{equation*}
	where $u^+,u^-$ and $\nu$ are as in definition \ref{approx-jump-set}.
\end{theorem}
For more on $BV$ functions one can see \cite{Ambrosio2}.

\noindent\textbf{Some preliminary results in hyperbolic system of conservation laws:}

 We denote $A(U)=DF(U)$ is the $n\times n$ Jacobian matrix of partial derivatives of $F$ where $F$ is the flux function for the system (\ref{HS1}). 
\begin{definition}\lab{SHS}
	We say that (\ref{HS1}) is strictly hyperbolic if for each $U\in\Omega$ the matrix $A(U)$ has $n$ real eigenvalues $\la_1(U)<\cdots<\la_n(U)$.
\end{definition}
By $r_{i}(U)$ we denote the right eigenvector corresponding to $i$-th eigenvalue.  
\begin{definition}
	We say $i$-th characteristic  field is genuinely non-linear if 
	\begin{equation}
	r_i\bullet \la_i(U)=D\la_i(U)\cdot r_i(U)> 0\ \mbox{ for all}\ U\in\Omega.
	\end{equation}
\end{definition}
\begin{definition}
	We say $i$-th characteristic  field is linearly degenerate if 
	\begin{equation}
	r_i\bullet \la_i(U)=D\la_i(U)\cdot r_i(U)= 0\ \mbox{ for all}\ U\in\Omega.
	\end{equation}
\end{definition}
\begin{definition}
	Let $U_0\in\Omega$. Suppose  $R_i(\si)(U_0)$ is the solution of the following ODE
	\begin{eqnarray}
	\frac{d}{d\si}R_{i}(\si)(U_0)&=&r_i(R_i(\si)(U_0)),\\
	R_i(0)&=&U_0.
	\end{eqnarray}
	We call the curve $R_i$ as $i$-rarefaction curve through $U_0$.	
\end{definition}
Next Theorem is due to Lax \cite{Lax}.
\begin{theorem}[Lax,\cite{Lax}]
	Suppose the system (\ref{HS1}) is strictly hyperbolic. Then, for each $U_0\in \Omega$ there exists a $\si_1>0$ and $n$-smooth curves $S_i:[-\si_1,\si_1]\rr\Omega$ together with scalar functions $\bar{\la}_i[-\si_1,\si_1]\rr\re$ such that
	\begin{equation}\lab{rh}
	F(S_i(\si))-F(U_0)=\bar{\la}_i(\si)\left(S_i(\si)-U_0\right)\ \mbox{ for }\si\in[-\si_1,\si_1].
	\end{equation} 
\end{theorem} 
Fix a $\si\in[-\si_1,\si_1]$. Suppose $U_-=U_0,U_+=S_i(\si),\la=\la_i(\si)$. Consider the Riemann data 
\begin{equation}\lab{rd}
U_0(x)=\left\{\begin{array}{lll}
U_- &\mbox{ for }&x<0,\\
U_+ &\mbox{ for }&x>0.
\end{array}\right.
\end{equation} 
Next, we consider the function $U:\re\times\re^+\rr\re^n$ defined as
\begin{equation}
U(x,t)=\left\{\begin{array}{lll}
U_0 &\mbox{ for }&x<t\la_i(\si),\\
S_i(\si) &\mbox{ for }&x>t\la_i(\si).
\end{array}\right.
\end{equation} 
Note that $U(x,t)$ satisfies the Rankine-Hugoniot equation (\ref{rh}), hence it is a weak solution to (\ref{HS1}) for the initial data (\ref{rd}). 
\begin{definition}[\textbf{Lax condition}]
	A weak solution $u$ is admissible if at every point $(x,t)$ of approximate jump, the left and right states $U_-,U_+$ and the speed $\la=\la_i(U_-,U_+)$ of the jump satisfy
	\begin{equation}
	\la_i(U_+)\leq\la\leq\la_i(U_-). 
	\end{equation} 	
\end{definition}
Next result is due to Foy \cite{Foy}.
\begin{theorem}[Foy, \cite{Foy}]
	If $i$-th characteristic field is genuinely non-linear then there is a $\si_0>0$ such that the solution (\ref{rd}) satisfies entropy inequality as well as the Lax admissible condition for all $\si\in[-\si_0,0]$.
\end{theorem}
For consistency of notation we denote $S_i(\cdot)$ by $S_i(\cdot)(U_0)$ and call it as \textit{$i$-th shock curve}. If $i$-th characteristic field is genuinely non-linear we choose the eigenvectors $r_i(U)$ so that $r_i\bullet\la_i\equiv1$. Further we choose the parametrization of the $i$-th shock and $i$-rarefaction curve such that 
\begin{eqnarray}
\frac{d}{d\si}\la_i(S_i(\si)(U_0))\equiv1,\ \frac{d}{d\si}\la_i(R_i(\si)(U_0))\equiv1,\\
\la_i(S_i(\si)(U_0))=\la_i(R_i(\si)(U_0))=\la_i(U_0)+\si.\lab{derivative}
\end{eqnarray}
Later on we utilize the following advantage of this parametrization,
\begin{equation}
U_0=S_i(-\si)\circ S_i(\si)(U_0)\ \mbox{ for all }\si, U_0.
\end{equation}
Equivalently, for $\si>0$ we have $U_0=S_i(-\si)(R_i(\si)(U_0))$.

\section{Main results}\lab{main-results}
\begin{theorem}\lab{theorem1}
	Let $f$ be a $C^2$ flux and $d\geq1$. Then there exists $u\in L^{\f}$, an entropy solution to (\ref{I1}) such that 
	\begin{equation}
	J\subsetneqq \mbox{cl}(J),
	\end{equation}
	where $J$ is defined as in (\ref{jump}) and $cl(J)$ denotes the closure of $J$.
\end{theorem}

\begin{theorem}\lab{theorem2}
	Let $F$ be a $C^2$ flux and $\de_0,T>0$. Let the system be strictly hyperbolic with $i$-th characteristic field is either genuinely non-linear or linearly degenerate for some $i\in\{1,\cdots,n\}$. Then there exists an entropy solution $U\in L^{\f}$ such that following holds
	\begin{enumerate}
		\item $	TV(U(\cdot,t))\leq\de_0$ for $t\in[0,T]$,
		\item $U$ is discontinuous on a set $\mathcal{B}$ and $\mathcal{B}\subsetneq cl(\mathcal{B})$.
	\end{enumerate}
Here $cl(J)$ denotes the closure of $J$.
  
\end{theorem}
\begin{remark}
 The solutions obtained in	Theorem \ref{theorem1} and Theorem \ref{theorem2} are discontinuous on a dense set.
\end{remark}
\section{Outline of the constructions}\lab{method}
\subsection{Outline of the construction for the scalar case}For the scalar case we first consider a strictly convex flux in one dimension and construct an entropy solution $u$ to (\ref{I1}) such that it is discontinuous over a dense set.
\begin{itemize}
	\item Employing the method of backward construction \cite{shyamcontrol,shyamop} we get a solution such that $y(x,t)$ has jump discontinuity over a dense subset of $[0,1]$. Here the function $y(x,t)$ comes from the Lax-Oleinik formula for strictly convex flux. 
	\item Then we analyze the one-sided limits of the entropy solution on the set $\mathcal{A}$ defined as follows.
	\begin{equation}
	\mathcal{A}:=\{\la(r_k,t)+(1-\la)(\rho(r_k),0);\la\in[0,1],k\in\mathbb{N}\},
	\end{equation}
	where $\{r_k;k\in\mathbb{N}\}$ is one enumeration of dyadic rational numbers in $[0,1]$ and the function $\rho$ is as in (\ref{rho}). In the first step, we construct the solution as a limit of a sequence of entropy solutions $u_n$. If $\mathcal{A}_n$ is the discontinuity set of $u_n$ then it follows that $\mathcal{A}\cap\mathcal{A}_n\subset \mathcal{A}\cap\mathcal{A}_{n+1}$. This leads to the conclusion that $u$ has discontinuity over $\mathcal{A}$.

\end{itemize} 
For a general $C^2$ flux in multi dimension we restrict ourselves in a particular direction $\xi$ and utilize the previous construction for $f\cdot\xi$. Then we extend the solution in multi-dimension by making it constant in any other orthogonal directions.
\subsection{Outline of the construction for the strictly hyperbolic  system}
Construction for system demands a different strategy than the one we have done for the scalar case since the backward construction for Lax-Oleinik solution is not available here. We proceed in the following way.
\begin{itemize}
	\item First we fix a characteristic field $\la_i$ which is either non-degenerate or linearly degenerate. For a fixed $U_0$ in $\Omega$ we consider 
	\begin{equation*}
	\Psi_i(\tau):=\left\{\begin{array}{lll}
	S_i(\tau) &\mbox{for}&\tau\in[-\tau_0,0],\\
	R_i(\tau) &\mbox{for}&\tau\in[0,\tau_0],
	\end{array}\right. \mbox{ for some }\tau_0>0,
	\end{equation*} 
	where $S_i,R_i$ are the \textit{shock curve} and \textit{rarefaction curve} respectively corresponding to $i$-th characteristic field (see section \ref{prelim} for definitions). We first choose $\{\tau_{k}\}_{1\leq k\leq m}$ such that
	\begin{equation*}
	\sum\limits_{k=1}^{m-1}|\la_i(\Psi_i(\tau_{k+1}))-\la_i(\Psi_i(\tau_k))|\leq C,
	\end{equation*} where $C$ is independent of $m$. Then we construct a piecewise constant admissible solution $U_m$ such that it takes values only on $\{\Psi_i(\tau_k)\}_{1\leq k\leq m}$. Then $TV(\la_i(U_m))$ is uniformly bounded. Then by Helly's Theorem we pass the limit and get an entropy solution.
	\item Note that in the previous step we have the flexibility to choose any sequence $\{\tau_k\}_{k\geq1}$ subject to $\sum\limits_{k=1}^{\f}|\tau_{k+1}-\tau_{k}|<\f$. In our construction we choose $\{\tau_k\}$ according to (\ref{11}) and  (\ref{12}). 
	\item We can prove that the solution we constructed is discontinuous over a dense set in a similar way as we have done for the scalar case. 
\end{itemize}

\section{Construction for scalar case}\label{scalar}
This section is devoted to the proof of Theorem \ref{theorem1}. We first show the $1$-dimensional case and then extend it to multi-dimension.
\subsection{Convex case in 1-dimension}\lab{scalar-1d}

\noindent\textbf{Construction:} Suppose $\left\{r_k\right\}_{k\geq1}$ is the enumeration of dyadic rational numbers in $[0,1]$. Define $\rho:\re\to\re$ as follows
\begin{equation}\lab{rho}
 \rho(x)=\left\{\begin{array}{lll}
                 0&\mbox{ for }&x\leq0,\\
                 \sum\limits_{k=1}^{\f}\frac{1}{2^k}\chi_{(-\f,r_k]}(x)&\mbox{ for }&0<x\leq1,\\
                 1&\mbox{ for }&x>1,
                \end{array}\right.
\end{equation}
where this $\chi_A$ is the characteristic function of the measurable set $A$. Note that this $\rho$ function is left continuous and increasing. We know from the Lax-Oleinik theory \cite{Lax,Oleinik} that if $u\in L^\f$ is the solution of
\begin{equation}
 \frac{\pa}{\pa t}u+\frac{\pa}{\pa x}f(u)=0,
\end{equation}
then we have the explicit formula $u(x,t)=G\left(\frac{x-y(x,t)}{t}\right)$ for $t>0$. This theory also tells us about incresing property of $y(x,t)$. Now we want to construct a solution show that at time $T>0$ this $y(x,T)$ coincides with our previously prescribed $\rho$ function. In order to do so we take help of backward construction (\cite{shyamop}) with special discretization. Then by some elementary analysis we prove that the solution $u(.,t)$ has discontinuiy on a dense set for a.e. $0<t<T$.

\begin{description}
\descitem{ \textit{Step(1)}\textbf{(Backward construction)}:}{step1} In this step, we construct an initial data $u_0$ to achive the solution $u(x,T)=G\left(\frac{x-\rho(x)}{T}\right)$ at time $t=T$ via backward construction. The methods and analysis of backward construction has been done in \cite{shyamcontrol,shyamop} to get optimal and exact controlability in conservation laws. Here we just mention the special sequences those are needed to be taken in order to prove the step \descref{step2}{2}. For full analysis we refer reader to \cite{shyamcontrol},\cite{shyamop}.
 \par We will consider the decomposition of domain $[0,1]$ as $[0,1]=\bigcup\limits_{j=1}^{k}(x_j,x_{j+1}]$ with $0=x_1<\cdots<x_{k+1}=1$ and $\left\{x_j;j=1,\cdots k+1\right\}=\left\{r_j;j=1,\cdots k+1\right\}$. Suppose $y_j=\rho(x_j)$ for $j=1,\cdots, k+1$. Define the piece-wise constant approximation of $\rho$ as
 \begin{eqnarray}
  \rho_k(x)=\left\{\begin{array}{lll}
                    x&\mbox{ if }& x\notin[0,1],\\
                    \sum\limits_{j=1}^{k}y_j\chi_{[x_j,x_{j+1})}&\mbox{ if }& x\in[0,1].
                   \end{array}
\right.
 \end{eqnarray}
By the choice of $r_k$ we have $|\rho(x)-\rho_k(x)|\leq\frac{1}{2^n}$ where $n$ is the largest integer less than or equal to $\log_2(k)$. Then we can construct a initial data $u^0_k$ such that at time $t=T$ the solution will be $u_k(x,T)=G\left(\frac{x-\rho_k}{T}\right)$. Also one can show that
$TV(f^\p(u^0_k(.)))\leq\frac{C}{T}$ with the constant $C>0$ which is independent of $k$. Then by Helly's theorem, we can extract a subsequence $\left\{f^\p(u^0_{k_l})\right\}_{l\geq1}$ such that it converges to some function $h(.)$ a.e. $x\in \re$ and in $L^1_{loc}$. From the Lipschitz continuity of $G$ we can obtain $u^0_{k_l}\to u^0=G(h)$ a.e. and in $L^1_{loc}$ as $l\to\f$. From $L^1$ contraction principle we get that $u_n\to u$ a.e. and in $L^1_{loc}$. Note that for each $k$ we have
\begin{equation}
 \int\limits_{\re}\int\limits_{0}^{T}u_k\phi_t(x,t)+f(u_k)\phi_x(x,t)dxdt+\int\limits_{\re}u^0(x)\phi(x,0)dx=0
 \end{equation} for all $\phi\in C_c^\f(\re\times[0,T))$.
 \begin{equation}
 \int\limits_{\re}\int\limits_{0}^{T}|u_k-m|\psi_t(x,t)+sgn(u_k-m)(f(u_k)-f(m))\psi_x(x,t)dxdt\geq0
 \end{equation}
 for all $0\leq\psi\in C^\f_c(\re\times(0,T))$ and $m\in\re$. Now by passing the limit $l\to\f$ we can show that $u$ is solution to
 \begin{eqnarray*}
  \frac{\pa}{\pa t} u+\frac{\pa}{\pa x}f(u)&=&0\mbox{ in }\re\times(0,T),\\
  u(x,0)&=&u^0(x)\mbox{ for }x\in\re.
 \end{eqnarray*}
 Since $u_k(x,T)=G\left(\frac{x-\rho_k(x)}{T}\right)$ we get $u(x,T)=G\left(\frac{x-\rho(x)}{T}\right)$. This completes step \descref{step1}{1}.

	\begin{figure}[ht]
	\centering
	\def\svgwidth{0.8\textwidth}
	\begingroup
	\setlength{\unitlength}{\svgwidth}
	\begin{picture}(1.05,0.66)
	\put(0,0.1){\includegraphics[width=0.85\textwidth]{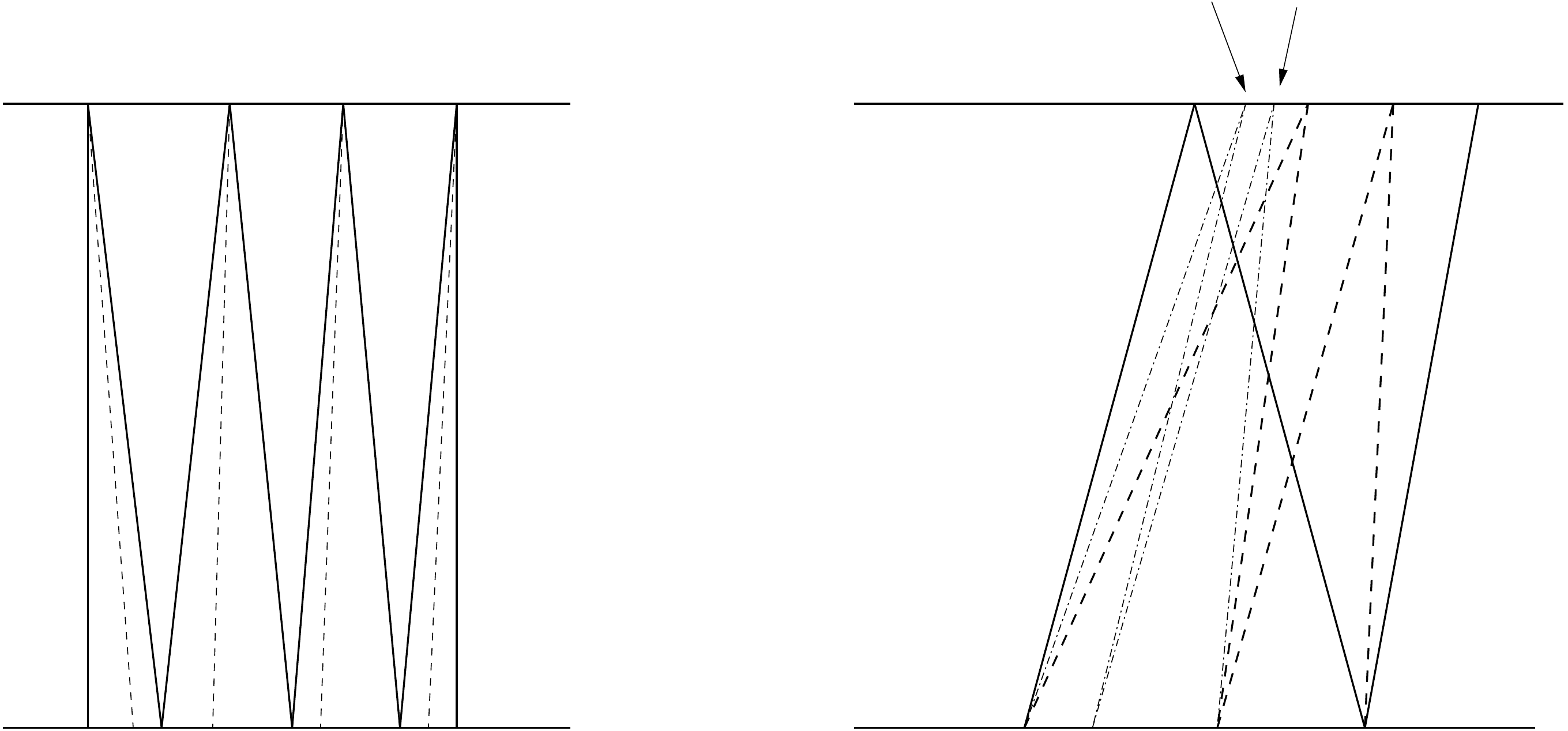}}%
	\put(0.15,0.54){\color[rgb]{0,0,0}\makebox(0,0)[lb]{\smash{$x_2$}}}%
	\put(0.05,0.54){\color[rgb]{0,0,0}\makebox(0,0)[lb]{\smash{$x_1$}}}%
	\put(0.22,0.54){\color[rgb]{0,0,0}\makebox(0,0)[lb]{\smash{$x_3$}}}%
	\put(0.3,0.54){\color[rgb]{0,0,0}\makebox(0,0)[lb]{\smash{$x_4$}}}%
	\put(0.1,0.08){\color[rgb]{0,0,0}\makebox(0,0)[lb]{\smash{$y_1$}}}%
\put(0.05,0.08){\color[rgb]{0,0,0}\makebox(0,0)[lb]{\smash{$y_0$}}}%
\put(0.18,0.08){\color[rgb]{0,0,0}\makebox(0,0)[lb]{\smash{$y_2$}}}%
\put(0.3,0.08){\color[rgb]{0,0,0}\makebox(0,0)[lb]{\smash{$y_4$}}}%
\put(0.26,0.08){\color[rgb]{0,0,0}\makebox(0,0)[lb]{\smash{$y_3$}}}
	\put(0.79,0.535){\color[rgb]{0,0,0}\makebox(0,0)[lb]{\smash{$x_{11}$}}}%
\put(0.885,0.535){\color[rgb]{0,0,0}\makebox(0,0)[lb]{\smash{$x_{22}$}}}%
\put(0.93,0.535){\color[rgb]{0,0,0}\makebox(0,0)[lb]{\smash{$x_{23}$}}}%
\put(0.98,0.535){\color[rgb]{0,0,0}\makebox(0,0)[lb]{\smash{$x_{12}$}}}%
\put(0.8,0.61){\color[rgb]{0,0,0}\makebox(0,0)[lb]{\smash{$x_{32}$}}}%
\put(0.87,0.61){\color[rgb]{0,0,0}\makebox(0,0)[lb]{\smash{$x_{33}$}}}%
\put(0.73,0.08){\color[rgb]{0,0,0}\makebox(0,0)[lb]{\smash{$y_{32}$}}}%
\put(0.68,0.08){\color[rgb]{0,0,0}\makebox(0,0)[lb]{\smash{$y_{11}$}}}%
\put(0.815,0.08){\color[rgb]{0,0,0}\makebox(0,0)[lb]{\smash{$y_{22}$}}}%
\put(0.92,0.08){\color[rgb]{0,0,0}\makebox(0,0)[lb]{\smash{$y_{12}$}}}%
\put(0.15,0.02){\color[rgb]{0,0,0}\makebox(0,0)[lb]{\smash{Fig. 1a}}}%
\put(0.8,0.02){\color[rgb]{0,0,0}\makebox(0,0)[lb]{\smash{Fig. 1b}}}%
	\end{picture}
	\endgroup
	\caption{ Fig. 1a shows the structure of $u_4$ obtained in backward construction. In this figure $y_i=\rho(x_i)$ for $1\leq i\leq4$. Fig. 1b shows a local behaviour of  $u_1,u_2,u_3$ near the line $L$ joining $(x_{11},T)$ and $(y_{11},0)$ where $u_1,u_2,u_3$ obtained in three consecutive steps of approximation.  Here $\{x_{mk}\}_{1\leq k\leq \Lambda_m}$ are the partitions in $m$-th step of the approximation for $1\leq m\leq 3$. Note that $x_{11}=x_{21}=x_{31}$.}
	\label{fig1}
\end{figure}
 \descitem{ \textit{Step(2):}}{step2} Now fix a time $0<t_0<T$ such that $u_{k_l}(x,t_0)\to u(x,t_0)$ a.e $x\in\re$ as $l\to\f$. Let $B$ be the subset of $[0,1]$ such that $|[0,1]\setminus B|=0$ and $u_{k_l}(y,t_0)\to u(y,t_0)$ for each $y\in B$. For each $j\geq1$ define
 \begin{equation*}
 z_j=x_j-\frac{(T-t_0)}{T}\left(x_j-y_j+\frac{1}{2^j}\right).
 \end{equation*}
 \begin{claim}\lab{cl1} For each fixed $j>2$ there exist two sequences $\left\{\bar{x}_m\right\}_{m\geq1},\left\{\tilde{x}_m\right\}_{m\geq1} $ such that
 \begin{eqnarray*}
  \bar{x}_{m}<z_j<\tilde{x}_m\mbox{ for }m\geq1,\mbox{ and }\lim\limits_{m\to\f}\bar{x}_m=z_j=\lim\limits_{m\to\f}\tilde{x}_m\\
  \liminf\limits_{m\to\f}f^\p(u(\tilde{x}_m,t_0))-f^\p(u(\bar{x}_m,t_0))\geq\frac{3}{2^{j+3}T}.
 \end{eqnarray*}
\end{claim}
\noindent\textit{Proof of the claim \ref{cl1}:} Fix a $\de>0$. Now choose $\bar{x},\tilde{x}$ from $B\setminus\bigcup\limits_{k=1}^{\f}\left\{z_k\right\}$ such that $z_j-\de<\bar{x}<z_j<\tilde{x}<z_j+\de$.
  Since $f^\p(u_{k_l}(y,t_0))\to f^\p(u(y,t_0))$ for $y=\bar{x},\tilde{x}$ there exists an $N_0$ such that for $l\geq N_0$ we have 
  \begin{equation}\lab{ineq1}
   |f^\p(u_{k_l}(y,t_0))- f^\p(u(y,t_0))|\leq\frac{1}{2^{j+2}T} \mbox{ for }y=\bar{x},\tilde{x}.
  \end{equation}
  Observe that there exists an $N_1$ such that if $k\geq N_1$ and $\{x_i\}_{i=1}^{k}$ the partition as mentioned in step(1) then $|x_i-x_{i+1}|\leq\frac{1}{2^{j+3}}$. Let $M=\max\{N_0,N_1\}$ and set $k_0=k_{N_0+1}$. Since $\bar{x},\tilde{x}$ are not equal to any of the $z_j$'s there exist $\bar{k}_1,\bar{k}_2,\tilde{k}_1,\tilde{k}_2$ such that 
 \begin{equation*}
  z_{\bar{k}_1}<\bar{x}<z_{\bar{k}_2}\mbox{ and }z_{\tilde{k}_1}<\tilde{x}<z_{\tilde{k}_2}
 \end{equation*}
 which implies 
 \begin{eqnarray*}
 \frac{x_{\bar{k}_1}-y_{\bar{k}_1}}{T} \leq f^\p(u_{k_0}(\bar{x},t_0))&\leq&\frac{x_{\bar{k}_2}-y_{\bar{k}_1}}{T},\\
 \frac{x_{\tilde{k}_1}-y_{\tilde{k}_1}}{T}\leq f^\p(u_{k_0}(\bar{x},t_0))&\leq&\frac{x_{\tilde{k}_2}-y_{\tilde{k}_1}}{T}.
 \end{eqnarray*}
 Therefore we have
  \begin{eqnarray*}
   f^\p(u_{k_0}(\bar{x},t_0))- f^\p(u_{k_0}(\tilde{x},t_0))&\geq& \frac{x_{\bar{k}_1}-y_{\bar{k}_1}}{T}-\frac{x_{\tilde{k}_1}-y_{\tilde{k}_1}}{T}\\
   &\geq&\frac{1}{T}(y_{\tilde{k}_1}-y_{\bar{k}_1})-\frac{1}{T}(x_{\bar{k}_1}-x_{\tilde{k}_2})\\
   &\geq&\frac{1}{2^{j}T}-\frac{1}{T}(x_{\bar{k}_1}-x_{\tilde{k}_2}).
 \end{eqnarray*}
 By choice of $k_0$ we have $|x_{\bar{k}_1}-x_{\tilde{k}_2}|\leq\frac{1}{2^{j+1}}$. This yields
 \begin{eqnarray*}
f^\p(u_{k_0}(\bar{x},t_0))- f^\p(u_{k_0}(\tilde{x},t_0))&\geq&  \frac{1}{2^{j}T}-\frac{1}{2^{j+3}T}\\
&\geq&\frac{7}{2^{j+3}T}>0.
 \end{eqnarray*}
Hence we get
\begin{equation}\lab{ineq2}
|f^\p(u_{k_0}(\bar{x},t_0))- f^\p(u_{k_0}(\tilde{x},t_0))|\geq\frac{7}{2^{j+3}T}.
\end{equation}
Finally inequalities (\ref{ineq1}) and (\ref{ineq2}) help us to conclude
\begin{eqnarray*}
 |f^\p(u(\bar{x},t_0))- f^\p(u(\tilde{x},t_0))|&=&|f^\p(u(\bar{x},t_0))- f^\p(u_{k_0}(\bar{x},t_0))+f^\p(u_{k_0}(\bar{x},t_0))- f^\p(u_{k_0}(\tilde{x},t_0))\\
 &+&f^\p(u_{k_0}(\tilde{x},t_0))- f^\p(u(\bar{x},t_0))|\\
&\geq&|f^\p(u_{k_0}(\bar{x},t_0))- f^\p(u_{k_0}(\tilde{x},t_0))|-|f^\p(u(\bar{x},t_0))- f^\p(u_{k_0}(\bar{x},t_0))|\\
&-&|f^\p(u_{k_0}(\tilde{x},t_0))- f^\p(u(\tilde{x},t_0))|\\
&\geq&\frac{7}{2^{j+3}T}-\frac{2}{2^{j+2}T}\\
&=&\frac{3}{2^{j+3}T}.
 \end{eqnarray*}
 By taking a sequence $\de_m\to0$ we get the sequences $\left\{\bar{x}_m\right\}_{m\geq1},\left\{\tilde{x}_m\right\}_{m\geq1}$ which  satisfies all the conditions prescribed in claim and we have
 \begin{equation*}
  |f^\p(u(\bar{x}_m,t_0))- f^\p(u(\tilde{x}_m,t_0))|\geq \frac{3}{2^{j+3}T}.
 \end{equation*}
 This yields 
 \begin{equation*}
    \liminf\limits_{m\to\f}|f^\p(u(\tilde{x}_m,t_0))-f^\p(u(\bar{x}_m,t_0))|\geq\frac{3}{2^{j+3}T}.
 \end{equation*}
 This claim shows that $f^\p(u(\cdot,t_0))$ is discontinuous at $z_j$ and so is $ u(\cdot,t_0)$. Hence discontinuity set of $u(\cdot,t_0)$ is a dense set.
  \end{description}
\subsection{General $C^2$ flux in multi dimension}\lab{scalar-multi}
 Suppose $f=(f_1,f_2,\cdots,f_d)\in C^2(\re,\re^d)$. Consider the set $X=\left\{f_1^{\p\p}=0\right\}$ and the following one dimensional scalar conservation laws
 \begin{eqnarray}\lab{1d_prob}
 \frac{\pa}{\pa s} v(\xi,s)+\frac{\pa}{\pa \xi}f_1(v(\xi,s))&=&0\mbox{ in }\re\times(0,\f),\\
  v(\xi,0)&=&v_0(\xi)\mbox{ for }\xi\in \re.
 \end{eqnarray} 
If interior of $X$ is non-empty then there exists an interval $I=(b_1,b_2)\subset X$ which means in the interval $I$ flux $f$ is linear. Hence for the initial data $v_0$ which takes value in $I$ and discontinuous at a dense set the entropy solution will have discontinuity at a dense set for all time $t>0$.
 Otherwise there exists an interval $I_1$ on which $f^{\p\p}\neq0$. Without loss of generality we may assume $f_1^{\p\p}>0$ on $J=\bar{I_1}$. This tells us $f|_J$ is uniformly convex. We have shown in previous section that there exists an initial data $v_0$ such that the entropy solution $v(.,t)$ is discontinuous on a dense set for $0<t<T$.
 Now set 
 \begin{equation}\lab{u0}
  u_0(x_1,x_2,\cdots,x_d)=v_0(x_1)\mbox{ for }x=(x_1,x_2,\cdots,x_d)\in\re^d.
 \end{equation}
Then the entropy solution $u$ to (\ref{I1}) corresponding to initial data $u_0$ as defined in (\ref{u0}) has discontinuity on a dense set.

This completes the proof of Theorem \ref{theorem1}.\qed
 

\section{Construction for strictly hyperbolic system}\lab{construction-system}
Let $\si_0>0$ be given. Consider the function $g:\re\rr\re$ defined as follows
 \begin{equation}
 g(x)=\left\{\begin{array}{lll}
 0 &\mbox{ for }x<0,\\
 \sum\limits_{k=1}^{\f}\frac{\si_0}{2^{k}}\chi_{(-\f,r_k]} &\mbox{ for }0<x<\si_0,\\
 \si_0 &\mbox{ for }\si_0<x,
  \end{array}\right.
  \end{equation}
  where $r_k$ defined in the section \ref{scalar}.
  Now fix a $U_0\in \Omega$. We first consider the case when $i$-th characteristic field is genuinely non-linear. We proceed as we have done for the scalar case. Here the main difference is that we do not have backward construction to build up a solution. But the method definitely motivates us to find a data and move forward with that.
  \begin{description}
  	\descitem{ \textit{Step 1:}}{step11} As we have seen in the construction for scalar case here  we first partition the line $\{(x,1);0\leq x\leq\si_0\}$ into $N$ parts $\{x_i^N;1\leq i\leq N\}$ such that  $|x_i^N-x_{i+1}^N|\leq\frac{\si_0}{2^k}$ for some $k\geq1$. Suppose $y_i^N=g(x_i^N)$ which partitions the line $\{(y,0);0\leq y\leq\si_0\}$. Set $\si_{j}^{N}=y_{j+1}^N-y_j^N$ and $\de_{j}^N=x_{j+1}^N-x_j^N$. Next consider following sequence of $U$.
\begin{eqnarray}
U_{m+\frac{1}{2}}&=&S_i\left(-\si_{m}^N+\sum\limits_{j=1}^{m-1}(-\si_{j}^N+\de_{j}^N)\right)\left(U_0\right),\lab{11}\\ 
U_{m+1}&=&S_i\left(\sum\limits_{j=1}^{m}(-\si_{j}^N+\de_{j}^N)\right)\left(U_0\right)\lab{12}.
\end{eqnarray}
It is easy to see that 
\begin{eqnarray}
\left|-\si_{m}^N+\sum\limits_{j=1}^{m-1}(-\si_{j}^N+\de_{j}^N)\right|&\leq&2\si_0,\\
\left|\sum\limits_{j=1}^{m-1}(-\si_{j}^N+\de_{j}^N)\right| &\leq&2\si_0.
\end{eqnarray}
Note that $\la_{i}(U_{m+\frac{1}{2}})<\la_i(U_m,U_{m+\frac{1}{2}})<\la_i(U_m)$. This yields there exists a $y_{m+\frac{1}{2}}^N\in(y_{m}^N,y_{m+1}^N)$ such that $x_m^N=y_{m+\frac{1}{2}}^N+\la_i(U_m,U_{m+\frac{1}{2}})$.

Next we define an initial data $V_0:\re\rr\re$ as follows
\begin{equation}
V_0(x)=\left\{\begin{array}{lll}
U_0 &\mbox{ if }&x<0,\\
U_{m}&\mbox{ if }& x_{m}^N<x<x_{m+\frac{1}{2}}^N,\\
U_{m+\frac{1}{2}} &\mbox{ if }& x_{m+\frac{1}{2}}^N<x<x_{m+1}^N,\\
U_0 &\mbox{ if }&\si_0<x.
\end{array}\right.
\end{equation}
Now note that 
\begin{eqnarray}
-\si_{m}^N+\sum\limits_{j=1}^{m-1}(-\si_{j}^N+\de_{j}^N)&=&-(y_{m+1}^N-y_1^N)+(x_{m}^N-x_1^N),\lab{21}\\
\sum\limits_{j=1}^{m-1}(-\si_{j}^N+\de_{j}^N)&=&-(y_{m}^N-y_1^N)+(x_{m}^N-x_1^N)\lab{22}.
\end{eqnarray}
Let $\{\bar{x}_{k}\}_{1\leq k\leq N_1}$ be another partition such that 
\begin{enumerate}
	\item $0=\bar{x}_1<\cdots<\bar{x}_{N_1}=\si_0$.
	
	\item $\{x_{k}\}_{1\leq k\leq N}\subset\{\bar{x}_{k}\}_{1\leq k\leq N_1}$.
\end{enumerate}
Suppose $\{\bar{U}_{m},\bar{U}_{m+\frac{1}{2}}\}_{1\leq m\leq N_1}$ are defined as in (\ref{11}),(\ref{12}). By the property we have $\bar{x}_1=x_1$ and $\bar{x}_{N_1}=x_N$. Suppose $x_{m_0}=\bar{x}_{m_1}$ for some $1\leq m_0\leq N$ and $1\leq m_1\leq N_1$. Then we have $y_{m_0}=\bar{y}_{m_1}$. From (\ref{11}),(\ref{12}),(\ref{21}) and (\ref{22}) we have $U_{m_0}=\bar{U}_{m_1}$. Next we observe that 
\begin{equation*}
\la_i(U_{m+\frac{1}{2}})=\la_i(U_m)-\si_m^N \mbox{ and }\la_i(U_{m+1})=\la_i(U_{m+\frac{1}{2}})+\de_m^N.
\end{equation*}
 This yields
 \begin{equation}
 TV(\la_i(V_0^N(\cdot)))\leq2\si_0.
 \end{equation} 
  By Helly's theorem, there exists a subsequence $\{V_0^{N_k}\}_{k\geq1}$ and $\la:\re\rr\re$ such that as $k\rr+\f$, $\la_i(V_0^{N_k}(x))\rr\la$ in $L^1_{loc}$ and a.e. $x\in\re$. Since $S_i([-\si_1,\si_1])$ and $\la_i\circ S_i([-\si_1,\si_1])$ are closed sets we can say $\la(\re)\subset \la_i\circ S_i([-\si_1,\si_1])$. Therefore, we can write $\la(x)=\la_i\circ S_i(\si(x))$ for some function $\si:\re\rr[-\si_1,\si_1]$. Thanks to the identity (\ref{derivative}) we have $\si(\cdot)\in BV(\re)$ and $TV(\si(\cdot)(\re)\leq C\si_0$ for some constant $C>0$ independent of $\si_0$. Since $S_i(\cdot)(U_0)$ is a smooth function we conclude that $TV(S_i(\si(\cdot)))\leq C_1\si_0$. Now we can choose $\si_0$ small enough so that $C_1\si_0<\de_0$. 
   Furthermore, we have for a.e. $(x,t)\in\re\times[0,1]$, $U^{N_k}(x,t)\rr U(x,t)$ as $k\rr+\f$.

 \descitem{ \textit{Step 2:}}{step12} Now fix a time $0<t_0<1$ such that $U^{N_k}(x,t_0)\to U(x,t_0)$ a.e $x\in\re$ as $k\to+\f$. Let $B_1$ be the subset of $[0,\si_0]$ such that $\mathcal{L}^1([0,\si_0]\setminus B_1)=0$ and $U^{N_k}(y,t_0)\to U(y,t_0)$ for each $y\in B_1$. For each $j\geq1$ define
\begin{equation*}
z_j=x_j-(1-t_0)\left(x_j-y_j+\frac{\si_0}{2^j}\right).
\end{equation*}
\begin{claim}\lab{cl2} For each fixed $j>2$ there exist two sequences $\left\{\bar{x}_m\right\}_{m\geq1},\left\{\tilde{x}_m\right\}_{m\geq1} $ such that
	\begin{eqnarray*}
	\dip	\bar{x}_{m}<z_j<\tilde{x}_m\mbox{ for }m\geq1,\mbox{ and }\lim\limits_{m\to\f}\bar{x}_m=z_j=\lim\limits_{m\to\f}\tilde{x}_m \ \mbox{with}\\
		\liminf\limits_{m\to\f}\la_i(U(\tilde{x}_m,t_0))-\la_i(U(\bar{x}_m,t_0))\geq\frac{3\si_0}{2^{j+3}}.
	\end{eqnarray*}
\end{claim}
This clam can be proved in same way we have done in the claim \ref{cl1} of step 2 for the scalar case. From the claim \ref{cl2} we conclude that the solution $U$ has discontinuity at each $(z_j,t_0)$ which is a dense subset of $(0,\si_1)\times\{t_0\}$. Hence the solution is discontinuous on a dense set in $[0,\si_0]\times[0,1]$ that is there is a set $A\subset[0,\si_0]\times[0,1]$ such that $U$ is discontinuous on $A$ and $\mathcal{L}^2\left(\left([0,\si_0]\times[0,1]\right)\setminus cl(A)\right)=0$. 
\end{description} 
Now we consider the case when characteristic field is linearly degenerate. Let us fix again a $U_0\in \Omega$. For this case we get the following identity (see \cite{Bressan3} for more details).
\begin{equation}\lab{ld1}
\la_i(U_0)=\la_i(R_i(\si)(U_0))=\la_i(S_i(\si)(U_0))\ \mbox{ for all }\si.
\end{equation} In this case we directly define the initial data. Let $\de>0$ be given. Let $\{\si_k\}$ be a sequence of real number such that the following holds,
\begin{equation}\lab{Es1}
|S(\si_k)(U_0)-U_0|\leq\frac{\de_0}{2^{k}}\ \mbox{ for all }k\geq1.
\end{equation} 
Next we define 
\begin{equation}\label{data}
V_0^N(x):=\left\{\begin{array}{lll}
U_0&\mbox{if}&x\leq0,\\
U_0+\sum\limits_{k=1}^{N}(S(\si_k)(U_0)-S(\si_{k-1})(U_0))\chi_{(-\f,r_k]}&\mbox{if}&0\leq x\leq 1,\\
U_0&\mbox{if}& x\geq 1.
\end{array}\right.
\end{equation}
Let $V^N(x,t)$ be the Lax entropy solution to (\ref{HS1}) for the data $V^N_0$. Then it follows from (\ref{ld1}) and (\ref{data}) that $V^N(x,t)=V_0^N(x-\la_i(U_0)t)$. From the estimate (\ref{Es1}) we can pass to the limit as $N\rightarrow \f$ and get the following initial data.
 \begin{equation}
 V_0(x)=\left\{\begin{array}{lll}
 U_0&\mbox{if}&x\leq0,\\
 U_0+\sum\limits_{k=1}^{\f}(S(\si_k)(U_0)-S(\si_{k-1})(U_0))\chi_{(-\f,r_k]}&\mbox{if}&0\leq x\leq 1,\\
 U_0&\mbox{if}& x\geq 1.
 \end{array}\right.
 \end{equation}
 and $TV(V_0)\leq\de_0$. Now by a similar argument as we have done for the genuinely non-linear case we get an entropy solution $V$ such that it has discontinuity on a dense set in $\re\times\re_+$ and $TV(V(\cdot,t))\leq\de_0$ for all $t>0$.
This completes the proof of Theorem \ref{theorem2}.\qed\\

\noindent\textbf{Acknowledgement.}  The first author would like to thank Inspire faculty-research grant\\ DST/INSPIRE/04/2016/000237.\\

\end{document}